\documentclass{article}
\usepackage{amsmath,amsfonts,amsthm}

\newtheorem{Problem}{Problem}
\newtheorem*{Algorithm}{Algorithm}
\newtheorem{Lemma}{Lemma}

\begin{document}

\title{Sequences of density $\zeta(k)-1$}
\author{William J. Keith}

%% Department of Mathematics, Drexel University, Philadelphia PA 19104
%% wjk26@drexel.edu

\maketitle

\section{INTRODUCTION}

At a recent post-seminar gathering, Herb Wilf casually mentioned to those of us assembled the fact that the quantities $\zeta(k)-1$ sum to 1, i.e., $$ \sum_{k=2}^{\infty} \left( \left( \sum_{i=1}^{\infty} \frac{1}{i^k} \right) - 1 \right) = 1 \, \text{.}$$

He then declared that when a sequence of nonzero positive numbers sums to 1, the entries of the sequence should well be interpretable as the probabilities of \emph{something}, and asked what it might be in this case.  The challenge here is to provide events as interesting as the numbers $\zeta(k)-1$ themselves: rich with internal relations and easy to describe.  One criterion for a really good answer to this challenge would be that the response supplies us with an intuitive feel for the relations among the quantities $\zeta(k)-1$.

The context of the discussion was number-theoretic, so after some consideration I formalized the challenge thus:

\begin{Problem} Partition $\mathbb{N}$ into sets $\{ {\cal{A}}_2, {\cal{A}}_3, \dots \}$ such that the asymptotic densities $\lim_{n \rightarrow \infty} \frac{1}{n} \# \left( {\cal{A}}_k \cap \{1,2,\dots , n \} \right) = \zeta(k)-1$.
\end{Problem}

This note provides one such partition.  There are, it will be seen, others, but {\ae}sthetics guides the choice presented here: the partition can be described in terms of components that build on each other, and we can produce an algorithm that can quickly determine the ${\cal{A}}_k$ to which any number belongs.

\section{THE CLASSES ${\cal{B}}_k$}

The key both to constructing the classes ${\cal{A}}_k$ and to proving that they have the required densities is to break down the $\zeta(k)-1$ thus:

\begin{equation}\label{ProbMatrix}
\begin{matrix}
\zeta(2) - 1 &= \frac{1}{2^2} &+ \frac{1}{3^2} &+ \frac{1}{4^2} &+ \frac{1}{5^2} &+ \dots \\
\zeta(3) - 1 &= \frac{1}{2^3} &+ \frac{1}{3^3} &+ \frac{1}{4^3} &+ \frac{1}{5^3} &+ \dots \\
\zeta(4) - 1 &= \frac{1}{2^4} &+ \frac{1}{3^4} &+ \frac{1}{4^4} &+ \frac{1}{5^4} &+ \dots \\
\zeta(5) - 1 &= \frac{1}{2^5} &+ \frac{1}{3^5} &+ \frac{1}{4^5} &+ \frac{1}{5^5} &+ \dots \\
\vdots          &          \vdots     &         \vdots       &      \vdots          &       \vdots        &    \\
          1        &=   \frac{1}{2}   &+    \frac{1}{6}  &+   \frac{1}{12}  & + \frac{1}{20} &+ \dots
\end{matrix}
\end{equation}

\noindent where the entries of the last line are the summations of the (infinite) columns above.  The strategy is now to construct classes each of which have asymptotic densities given by the individual entries of the table.  As those entries are all unit fractions, our task is much simplified: ask for a set of integers of some peculiar irrational density and intuition would not necessarily leap to service; but ask for one-sixth of the integers, and any residue class mod 6 will do.  So we begin in this section with that last row, by producing classes ${\cal{B}}_j$ of densities $\frac{1}{2}, \frac{1}{6}, \frac{1}{12}, \dots$, representing each column.  In the next section we break each ${\cal{B}}_j$ apart into sections ${\cal{B}}_{j,k}$ of densities $\frac{1}{j^k}$, $2 \leq j,k < \infty $, to produce the individual entries of the table.  The finishing stroke is then to collect across rows to form the ${\cal{A}}_k$.

The constructions we use are often explained nicely by recasting them in terms of the \emph{factorial-base} representation of a number, its expansion as $x = a_1 \cdot 1! + a_2 \cdot 2! + \dots + a_n n!$, $0 \leq a_i \leq i$).  This expansion is unique, as the student unfamiliar with this base system is encouraged to try proving on their own.  For more on such expansions, see \cite{Barwell}.

The class ${\cal{B}}_2$ must have asymptotic density $\frac{1}{2}$.  Let us select it to consist of those $x \equiv 1$ mod 2.

Modulo 6, the classes 1, 3, and 5 have been assigned, so for ${\cal{B}}_3$ we select the residue class 2 mod 6.

Modulo 24, we require two residue classes in order that ${\cal{B}}_4$ should have density $\frac{1}{12}$.  We select the first two available: those $x \equiv 4$ or 6 mod 24.

Modulo 120, we require six residue classes in order that ${\cal{B}}_5$ should have density $\frac{1}{20}$, and the first six available are 10, 12, 16, 18, 22, and 24.  

Notice that we cannot use modulus 20 for ${\cal{B}}_5$: every residue class mod 20 has already had some subprogressions assigned.  Instead, we will use the moduli $j!$.  Since for $m<j$, residue classes modulo $m!$ are the same as $\frac{j!}{m!}$ equally spaced residue classes modulo $j!$, all previous assignments can be viewed as multiple assignments modulo $j!$.  For each ${\cal{B}}_j$, we select the first $(j-2)!$ residue classes that have not yet been assigned.  What are these in general?

\begin{Lemma}\label{BAssignment} The assigned residue classes for ${\cal{B}}_m$ are $$\{ x \equiv \, 0! + \sum_{j=1}^{m-2} a_j \cdot j! \, \text{mod} \, m! \, , 1 \leq a_j \leq j \} \, \text{.}$$
\end{Lemma}

\begin{proof} In the examples above, by the time we have selected component residue classes for ${\cal{B}}_j$ we have assigned membership in some ${\cal{B}}_i$ to the natural numbers $1, 2, \dots , (j-1)!$, as part of residue classes either modulo $j!$ or modulo divisors thereof (previous factorials).  Let us proceed by induction, and suppose that this has been the case through ${\cal{B}}_{m-1}$: for each $b \leq m-1$, the numbers $1, 2, \dots (b-1)!$ have been assigned to some ${\cal{B}}_j$, $j \leq b$, as part of residue classes modulo $j!$.

Given this induction hypothesis, what is the first residue class modulo $m!$ available for assignment to ${\cal{B}}_m$?  None of 1 through $(m-2)!$, to begin with: each of those classes were assigned to ${\cal{B}}_{m-1}$ or some previous class.  We do know that none of the next $(m-2)!$ residue classes mod $m!$, $(m-2)!+1$ through $(m-2)!+(m-2)!$, are assigned to ${\cal{B}}_{m-1}$.  Examine this segment.

Since $(m-2)! \equiv 0$ mod $(m-2)!$, our assumption tells us that the numbers $(m-2)!+1$ through $(m-2)!+(m-3)!$ must have been assigned to ${\cal{B}}_{m-2}$ or earlier classes, as residues modulo $(m-2)!$ or divisors thereof.  However, none of $(m-2)!+(m-3)!+1$ through $(m-2)!+(m-3)!+(m-3)!$ were assigned to ${\cal{B}}_{m-2}$, and this segment is wholly contained within those numbers we already knew were not assigned to ${\cal{B}}_{m-1}$.  At each stage of the analysis, we have an unbroken segment of assigned classes for the numbers $1, 2, \dots, (m-2)! + (m-3)! + \dots + (m-i)!$, and we know that the next $(m-i)!$ numbers are not assigned to ${\cal{B}}_i$ through ${\cal{B}}_{m-1}$.  At the last step, we find that the last number assigned is $(m-2)! + (m-3)! + \dots + 2! + 1!$, an odd number in class ${\cal{B}}_2$.

The very next number is in the segment we know has not been assigned to ${\cal{B}}_2$ or any previous class, so it is available.  Thus the first available residue class modulo $m!$ to assign to ${\cal{B}}_m$ is $(m-2)! + (m-3)! + \dots + 2! + 1! + 0!$, a construction for which there are unfortunately multiple labels and notations in current usage.  The relevant sequence ${1, 2, 4, 10, 34, 154, 874, \dots }$ is Sloane's A003422 (\cite{Sloane}) missing a leading zero term, and is there called the \emph{left factorial} and denoted $!n = 0! + \dots +(n-1)!$.  We will use this notation, so that the first unassigned class modulo $m!$ is $!(m-1)$.

Now let us determine the remainder of the assigned classes mod $m!$.  At the next-to-last step of our analysis previously, assuming $m \geq 4$ so that such a step exists, we determined that $!(m-1) - 2$ must have been an assigned class modulo $3!$, that is, $!(m-1)-2 \equiv 2$ mod 6.  Then $!(m-1)$ is not such a class and neither is $!(m-1)+2$: these are the two even classes mod 6 that were not assigned to ${\cal{B}}_3$.  Nor can they be assigned to ${\cal{B}}_i$ for $i \geq 3$ as part of residue classes modulo multiples of $3!$.  Thus, we can select these two as our first two residue classes modulo $m!$ for ${\cal{B}}_m$.  If $m \geq 5$, we can ascend backward another step: there we find that $!(m-1)-4$ and $!(m-1)-6$ had to be the assigned residue classes mod $4!$.  We already knew $!(m-1)$ and $!(m-1)+2$ were not in those classes, and here we find that they can be joined by $!(m-1)+6$ or $+12$, and $!(m-1)+2+6$ or $+12$.  These are the six residue classes modulo 24 that were left after we assigned two of them for ${\cal{B}}_4$.

At the $i$-th step of the ascent, we find that to all the entries we previously picked, we can add as many as $i-1$ multiples of $i!$.  Thinking about it in the factorial-base representation, using $a_1 a_2 a_3 \dots$ to denote $a_1 \cdot  1! + a_2 \cdot 2! + a_3 \cdot 3! + \dots$, the assigned classes will be $02111\dots + 0 b_2 b_3 b_4 b_5 \dots$, with $0 \leq b_i \leq i-1$.  And these are exactly the numbers claimed in the Lemma.
\end{proof}

\subsection{A Fast Membership Algorithm}

Having assigned the whole numbers to their various ${\cal{B}}_m$, we would like a means of determining which ${\cal{B}}_m$ a given $x$ belongs to without constructing every ${\cal{B}}_m$ until we assign $x$.  As it turns out, systematically selecting the first available classes at each step helps us write a computationally simple algorithm to do so.

Write $x$ in the factorial-base representation, $x = x_1 \cdot 1! + x_2 \cdot 2! + x_3 \cdot 3! + \dots + x_j \cdot j!$, $0 \leq x_i \leq i$.  We know that if the least positive residue of $x$ modulo any $n!$ is less than $(n-1)!$, $x$ must have been assigned to one of the ${\cal{B}}_{j \leq n}$, because from our proof above, the assigned residue classes for the ${\cal{B}}_{j \leq n}$ were an unbroken string from 1 to $(n-1)!$.

A "carry" in factorial-base addition occurs when the sum of the coefficients on the two summands on $i!$ is at least $i+1$, since $(i+1)i! = 1 \cdot (i+1)!$.  But the Lemma tells us that the residues assigned to ${\cal{B}}_m$ are of the form $02111\dots + 0 b_2 b_3 b_4 b_5 \dots$, with $0 \leq b_i \leq i-1$.  So in order for a 0 to ever occur at the $i!$ place of the factorial-base expansion of a residue assigned to ${\cal{B}}_m$, when $1 < i \leq m-1$, we need a very specific summand.  We need $b_2 = 1$, so that we carry from $2!$, leaving a 0.  We then need $b_3 = 2$, so that the sum of coefficients on $3!$ is 1 from $!(m-1)$, plus 2 from $b_3$, plus a carried 1, making 4, so we carry 1 and leave a 0 again.  This string has to continue up to the $i!$ term where we desire a 0.

Thus we have one of two cases that will diagnose where $!(m-1)$ ends for a given $x$.  We have that the least positive residue of $x$ modulo $m!$ is $\leq (m-1)!$.  Possibly $x \equiv (m-1)!$ modulo $m!$ exactly, so the string of $x_i$ starts with a possibly empty string of 0s that terminates with a 1 in position $x_{m-1}$: we added the largest possible value for every $b_i$, and carried at every step of the addition.  If we did not carry, then there is some smallest $m$ such that $x_{m-1} = 0$ but no previous $x_i = 0$, except for possibly an initial string that does \emph{not} terminate with a 1.  Write out $x$ in the factorial-base representation, check for the first behavior, and if it doesn't happen find the first 0.  The index where whichever of these occurs, occurs, is $m-1$.

\section{THE CLASSES ${\cal{A}}_k$}

Now that we have identified the classes ${\cal{B}}_m$, there remain the tasks of breaking them up into subclasses to give us the individual entries of Table \ref{ProbMatrix}, and collecting the corresponding parts "horizontally" to form the ${\cal{A}}_k$.

Each of the columns' individual entries decrease in geometric progression.  The class ${\cal{B}}_2$ has density $\frac{1}{2} = \frac{1}{2} \left( \frac{1}{2} + \frac{1}{4} + \frac{1}{8} + \dots \right) $.  The class ${\cal{B}}_3$ has density $\frac{1}{6} = \frac{1}{6} \left( \frac{2}{3} + \frac{2}{9} + \frac{2}{27} + \dots \right) $.  The class ${\cal{B}}_4$ has density $\frac{1}{12} = \frac{1}{12} \left( \frac{3}{4} + \frac{3}{16} + \dots \right) $, etc.

These geometric series give us the subclasses we need.  We assign the earliest fraction $\frac{m-1}{m}$ of ${\cal{B}}_m$ to the first subset, the next fraction $\frac{1}{m} \frac{(m-1)}{m}$ to the next set, and so forth.  This can be done by a one-step digit test.  For an example from familiar territory, to obtain $\frac{9}{10}$ of the integers, take those ending in 1 through 9.  Then to obtain $\frac{9}{100}$ of the integers by choosing $\frac{9}{10}$ of the remaining tenth, take those that end in 0 but have one of 1 through 9 in the tens place.

Since ${\cal{B}}_m$ is made up of residue classes modulo $m!$, rewrite $x$ as $x = x_0 + m! \cdot \left( b_0 m^0 + b_1 m^1 + b_2 m^2 + \dots \right) $, $0 \leq b_i \leq m-1$.  If $n$ is the smallest number such that $b_n \neq m-1$ assign $x$ to the set ${\cal{B}}_{m,n+2}$.  In table \ref{ProbMatrix}, $x$ is part of the set with asymptotic density given by the entry in column $m$, row $n+1$, which represents an entry in the sum for ${\cal{A}}_{n+2}$.

Thus, for example, we break up ${\cal{B}}_2$ as follows: ${\cal{B}}_{2,2} = \{ x \equiv 1 \, \text{mod} \, 4(=2 \cdot 2!) \}$, with density $\frac{1}{4}$; ${\cal{B}}_{2,3} = \{ x \equiv 3 \, \text{mod} \, 8(=4 \cdot 2!) \}$, with density $\frac{1}{8}$; ${\cal{B}}_{2,4} = \{ x \equiv 7 \, \text{mod} \, 16(=4 \cdot 2!) \}$, etc.  We break up ${\cal{B}}_3$ into fractions of size $\frac{2}{3}, \frac{2}{9}$, etc.: ${\cal{B}}_{3,2} = \{ x \equiv 2 \, \text{or} \, 8 \, \text{mod} \, 18(=3 \cdot 3!) \}$, ${\cal{B}}_{3,3} = \{ x \equiv 14 \, \text{or} \, 32 \, \text{mod} \, 54(=9 \cdot 3!) \}$, etc.

We now sum up by collecting these subsets across rows: ${\cal{A}}_k = {\bigcup \atop m} {\cal{B}}_{m,k}$.  Since the ${\cal{B}}_{m,k}$ are disjoint and the series of partial sums of their densities converges absolutely to  $\zeta(k)-1$, ${\cal{A}}_k$ has exactly the required asymptotic density.

The algorithm to determine which class ${\cal{A}}_k$ a given $x$ belongs to is thus:

\begin{Algorithm}
Write $x = a_1 \cdot 1! + a_2 \cdot 2! + a_3 \cdot 3! + \dots $, $0 \leq a_i \leq i$.

Let $j$ be the smallest index such that $a_j \neq 0$.  If $a_j=1$, then $m=j+1$.  If $a_j \neq 1$, then $m$ is the smallest number such that $a_i \neq 0$ for $j \leq i < m-1$, and $a_{m-1} = 0$.  Then $x \in {\cal{B}}_m$ .

Rewrite $x$ as $x = x_0 + m! \cdot (b_0 m^0 + b_1 m^1 + b_2 m^2+ \dots )$, $0 \leq b_i < m$, $0 < x_0 < m!$.  Let $n$ be the smallest index such that $b_n \neq m-1$.  Then $k = n+2$, and $x \in {\cal{A}}_k$.
\end{Algorithm}

\section{DATA AND SPECULATION}

The first few of the sets ${\cal{B}}_k$ and ${\cal{A}}_k$ are listed below.

${\cal{B}}_2 = \{ 1, 3, 5, 7, 9, 11, 13, 15, \dots \}$.

${\cal{B}}_3 = \{ 2, 8, 14, 20, 26, 32, \dots \}$.

${\cal{B}}_4 = \{ 4, 6, \quad 28, 30, \quad 52, 54, \quad 76, 78, \dots \}$.

${\cal{B}}_5 = \{ 10, 12, 16, 18, 22, 24, \qquad 130, 132, 136, 138, 142, 144, \dots \}$.

${\cal{B}}_6 = \{ 34, 36, 40, 42, 46, 48, 58, 60, 64, 66, 70, 72, \dots, 120, \dots \} $.

${\cal{A}}_2 = \{ 1,2,4,5,6,8,9,10,12,13,16,17,18,20,21,22,24,25,26, \dots \}$.

${\cal{A}}_3 = \{ 3,11,14,19,27,32,35,43,51,59,67,68,75,76,78,83,86, \dots \}$.

${\cal{A}}_4 = \{ 7,23,39,50,55,71,87,103,104,119,135,151,167,183,199, \dots \}$.

${\cal{A}}_5 = \{ 15,47,79,111,143,158,175,207,239,271,303,320,335, \dots \}$.

More terms, for the series up to ${\cal{A}}_{10}$, have been submitted to Sloane's database \cite{Sloane} and should be available by the time of publication of this note.

In terms of their component residue classes, the sets are:

${\cal{B}}_2 = \{ x \equiv 1 \, \text{mod} \, 2 \}$.

${\cal{B}}_3 = \{ x \equiv 2 \, \text{mod} \, 6 \}$.

${\cal{B}}_4 = \{ x \equiv 4, 6 \, \text{mod} \, 24 \}$.

${\cal{B}}_5 = \{ x \equiv 10, 12, 16, 18, 22, 24 \, \text{mod} \, 120 \}$.

${\cal{B}}_6 = \{ x \equiv 34, 36, 40, 42, 46, 48, 58, 60, 64, 66, 70, 72, \dots, 120 \, \text{mod} \, 720 \} $.

${\cal{A}}_2 = \{ x \equiv 1  \, \text{mod} \, 4; \, 2,8 \, \text{mod} \, 18; \, 4,6,28,30,52,54  \, \text{mod} \, 96 \dots \}$.

${\cal{A}}_3 = \{ x \equiv 3 \, \text{mod} \, 8; \, 14,32  \, \text{mod} \, 54; \, 76,78,172,174,268,270 \, \text{mod} \, 384 \dots \}$.

${\cal{A}}_4 = \{ x \equiv 7  \, \text{mod} \, 16; \, 50, 104  \, \text{mod} \, 162; \, 364, 366, \dots ,1134 \, \text{mod} \, 1536 \dots \}$.

Some variations of this construction could be explored.  Building ${\cal{B}}_j$, we had to make a series of choices.  We chose the odd numbers to form ${\cal{B}}_2$, but the even numbers would have worked just as well.  We chose the residue class $x \equiv 2$ mod 6 for the class ${\cal{B}}_3$, but could just as easily have chosen $x \equiv 4$ or 6.  Choosing the first available residue class for each factorial modulus gave us classes ${\cal{B}}_m$ that could be determined with a short algorithm, providing a tidy answer to the original problem.  However, other choices lead to answers with different features.

For example, suppose that at each step we choose the \emph{last} available residue class modulo $j!$ to construct ${\cal{B}}_j$.  Choose those $x \equiv 2$ mod 2 for ${\cal{B}}_2$.  Of the three available odd residue classes mod 6, choose $x \equiv 5$ mod 6 for ${\cal{B}}_3$.  Use the classes $x \equiv 19, 21$ mod 24 for ${\cal{B}}_4$.  With such a decision procedure, some numbers are never assigned to any ${\cal{B}}_j$ at all!  

The "missed set" is ${1, 3, 7, 9, 13, 15, 25, 27, 31, 33, \dots = 1 + \sum a_i i!, 0 \leq a_i < i}$.  All of the ${\cal{B}}_j$ will still have the same densities, summing to 1, so the missed set is "small" in that it is of asymptotic density 0.  In fact, given any arithmetic progression $X$ mod $Y$, some subprogression will be assigned to a ${\cal{B}}_j$.  (Show it!)  On the other hand, at any given point the missed set may be rather large for some purposes: $1/n$ of the numbers smaller than $n!$ are permanently unassigned.

These two assignment procedures are in some sense on opposite poles of an entire ensemble of possible procedures.  An interested reader might burrow a layer deeper than we have: assign some straightforward process for describing and choosing assignment procedures to construct each ${\cal{B}}_j$ from residue classes mod $j!$ and examine the resulting ensemble of all possible constructions.  Does it possess any striking structural features?

An early stab at this problem involved the powerfree numbers.  Squarefree numbers are those whose prime factors are all distinct: they have density $\frac{1}{\zeta(2)}$ in the whole numbers.  Cubefree numbers have factorizations in which no prime factor is repeated more than twice, so squarefree  numbers are also cubefree.  Cubefree numbers are of density $\frac{1}{\zeta(3)}$, and those which are cubefree but not squarefree are of density $\frac{1}{\zeta(3)} - \frac{1}{\zeta(2)}$.  Those that are 4th-powerfree but not cubefree are of density $\frac{1}{\zeta(4)} - \frac{1}{\zeta(3)}$, and so on.  Every whole number is $(k+1)$st-powerfree but not $k$th-powerfree for some $k$, the largest exponent in its prime factorization.

These densities seemed to suggest the possibility of defining a simple probability distribution on $\mathbb{N}$ that assigned a total probability of $\zeta(k) - 1$ to the event that a random integer variable would be $(k+1)$st-powerfree but not $k$th-powerfree.  It is easy to provide such a distribution by brute force -- say, giving the $n$th number which is $(k+1)$st-powerfree but not $k$th-powerfree a probability of $(\zeta(k) - 1) \cdot \frac{1}{2^n}$.\footnote{Typo in this line corrected in preprint at kind communication from Michael Lugo.}  But this does not seem to be a particularly illuminating illustration of any relations between the values $\zeta(k)$; indeed, such a construction works with any partitioning of the integers, and any set of probabilities instead of $\zeta(k)-1$.  A distribution based on the factorization of $x$ would seem much more natural; can a simple one be produced?

In closing, I would like to mention a personal recollection.  In general, the density of ${\cal{B}}_j$ is $\frac{1}{j(j-1)}$. Seeing it written on the wall of my office, I was reminded of the first place I ever encountered Leibniz' summation of the triangular series: William Dunham's delightful popular-mathematics text, \emph{Journey Through Genius} \cite{JTG}.  I read this book in high school, and it motivated in considerable part my decision to pursue mathematics in college.  The present note gives me an opportunity to thank Mr. Dunham sincerely for the service.

\end{document}